# Optimization of the closed-loop controller of a discontinuous capsule drive using a neural network

*Sandra Zarychta[1,2,\*], Marek Balcerzak[1], Volodymyr Denysenko[1], Andrzej Stefański[1], Artur Dąbrowski[1], Stefano Lenci[2]*

[1] *Division of Dynamics, Lodz University of Technology, Stefanowskiego 1/15, 90-924, Lodz, Poland*

[2] *Department of Civil and Building Engineering, and Architecture, Polytechnic University of Marche, via Brecce Bianche, 60131, Ancona, Italy*

## Abstract

In this paper, construction of a neural-network based, closed-loop control of a discontinuous capsule drive is analyzed. The foundation of the designed controller is an optimized open-loop control function. A neural network is used to determine the dependence between the open-loop controller's output and the system's state. Robustness of the neural controller with respect to variation of parameters of the controlled system is analyzed and compared with the original, optimized open-loop control. It is expected that the presented method can facilitate construction of closed-loop controllers of systems, for which other methods are not effective, such as non-smooth or discontinuous ones.

Keywords: *Optimal control, closed-loop, neural network, capsule drive, discontinuous system, non-smooth system*

## 1. Introduction

Capsule robots, often abbreviated as capsubots, are a class of micro capsule-shaped robots able to explore fields which are normally inaccessible to humans [1]. A particularly interesting subset of capsubots are devices propelled by an internal mechanical oscillator. The vibrating mass produces inertia forces which, in the presence of friction, allow to move the whole capsule. In this approach, external moving parts such as wheels, tracks, robotic legs or arms are no longer necessary [1, 2]. Such capsubots are very interesting from the practical point of view, due to their enormous potential in medicine, engineering and other areas [1, 2]. Moreover, their rich dynamics [2] which encompasses phenomena such as impacts, dry friction etc. is a broad research topic itself.

A lot of research concerning analysis and design of the capsubots has already been performed. For instance, the work [2] presents a mesoscale prototype of a self-propelled vibro-impact capsule system, as well as its optimization in terms of the average progression velocity, energy efficiency and power consumption. In the paper [3], development of control strategies for a capsubot is proposed. The inner mass of the device consists of two masses placed at both ends of the cylindrical rod is surrounded by a motor housing with a coil held in a shell [3]. Another interesting example is presented in [4], where the proposed vibro-impact capsule contains a harmonically excited internal oscillator impacting a massless plate suspended on a spring. Such arrangement causes that the resultant horizontal force acting on the capsule is not symmetric, which in the presence of dry friction enables the system to move forward. A more detailed study on the control function, where maximization of the rate of progression, as well as optimization of the energy consumption has been included, is shown in [5–7]. Moreover, the detailed bifurcation analysis of a vibro-impact system with use of path-following methods, accompanied by experimental investigation, can be found in [8]. The paper [9] describes a downscaled self-propelled vibro-impact capsule system with the ability to move precisely in a limited space, with a size equal to a market-leading gastrointestinal capsule endoscope. The capsule system includes two impact constraints with a linear bearing holding a T-shape magnet situated between them that restricts its linear motion. The dynamic analysis of the prototype is described, as well as optimization of the progression speed and minimization of the required propulsive force.

Apart from the vibro-impact systems, a different layout can be utilized in capsule drives. In particular, an interesting modification of the vibro-impact capsule can be obtained by replacing the mass-on-spring oscillator with a pendulum. In such case propulsion of the capsule is caused by interactions between friction, inertia forces produced by swinging motion of the pendulum and the contact force between the capsule and the underlying surface [10]. Such arrangement seems to make dynamics of the system somewhat more complex, because the contact force is dependent on oscillations of the pendulum. Periodic locomotion principles and nonlinear dynamics of a pendulum-driven capsule system are investigated in the works [10–12]. Moreover, the paper [12] also includes a motion generation strategy in the presence of visco-elasticity. Design and parameters optimization of a pre-designed control function profile for the pendulum capsule system is considered in [12–14].

Existing methods of controllers design applicable to capsubots use various approaches, including open-loop control, closed-loop feedback linearization or neural networks. The work [1] presents three control approaches for capsubot. The first one involves an open-loop control, whereas the second utilizes a closed-loop control with partial feedback linearization technique based on trajectory tracking. The last one, called simple switch control, is a combination of the previous methods. The control profiles learned from the open and closed-loop control are used to move the capsubot effectively in the desired direction. In [3], the authors present a strategy of controlling cylindrical rod composed of two stages. The first one assumes desired trajectory generation, whereas the second one focuses on inner mass closed-loop control for a given desired trajectory with a partial feedback linearization approach. An adaptive trajectory tracking control method for a vibro-driven capsule system is described in [13]. The implementation of an auxiliary input control variable, establishing the non-collocated feedback loop, is constructed to cope with the parametric uncertainties. A comparison of the proposed approach with the classical one has been performed with the use of a closed-loop feedback tracking control system. Improvements of this method have been shown in the paper [15]. The novel approach focuses on adding a neural network approximator and a robust compensator to an auxiliary control variable. The proposed design method with multi-layer neural networks and variable strategy structure, as well as an adaptive tracking control scheme, copes well with uncertainties such as parameters not known a priori, approximation errors, and disturbances [15]. A novel, Fourier series based method of open-loop optimal control estimation, applicable for discontinuous systems such as capsule drives, has been described in the work [16]. In the paper [17] a problem of crossing a circular fold by a capsule robot is concerned. The path following techniques have been utilized and the COCO software was used in numerical studies.

The paper [18] describes speed optimization of the self-propelled capsule robot [9] in the varying friction environment between the device and its supporting surface with the use of Six Sigma and Multi-Island genetic algorithms, where the Monte Carlo approach has been utilized for validation. Other examples of multi-objective optimization with Six Sigma as a controller for the genetic algorithm are presented in works [19, 20], whereas the reliability analyses with the use of a Monte Carlo algorithm is described in papers [21, 22].

It seems that there is little research on application of neural networks in optimal control of capsule drives or other similar systems. Possibly, direct application of Reinforcement Learning technique [23–25] could be used to obtain an approximation of the optimal closed-loop controller, but such approach would require a lot of time and computational resources [23–25]. However, a simpler option is proposed in this work. Providing that an optimal open-loop control is determined, a neural network approach can be used to determine the dependence of the controller's output and the corresponding states of the controlled system. In such a manner, a closed-loop controller, whose action reflects the way in which the open-loop optimal control works, can be obtained.

The aim of the study is to test and evaluate the aforementioned concept. For this purpose, an approximation of the open-loop optimal control of a pendulum capsule drive is performed by means of the method described in [16]. After that, a neural network is used to determine dependency between the output value of the optimized open-loop controller and the corresponding states of the capsule system. In such a manner, a closed-loop controller is obtained. Finally, performance and robustness of the closed-loop neural controller is compared with the original open-loop one. Results show that the neural controller maintains efficiency of the original and offers much better robustness against uncertainty of the friction coefficient of the controlled system, which actually is one of the main limitations in the use of open-loop controllers.

Authors believe that such solution can be an interesting option in design and optimization of controllers used in mechanical systems, including discontinuous ones. Moreover, it is expected that proposed method will facilitate construction of closed-loop controllers of systems, for which an optimal open-loop control is available.

## 2. Background

The subject of this research is the pendulum capsule drive. This chapter presents a brief description of the system along with approximation of its open-loop optimal control. Information presented below is a foundation for the new, neural network based, closed-loop controller of the device. A scheme of the pendulum capsule drive is presented in Fig. 1.

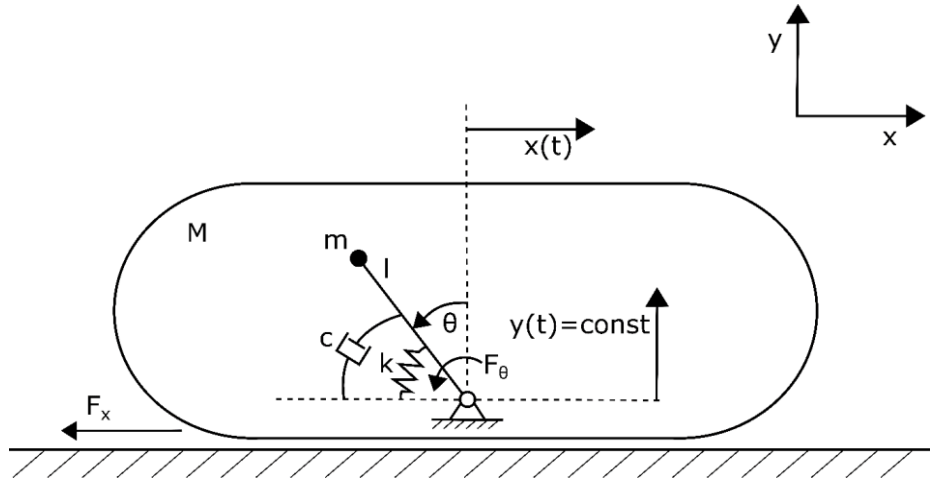

*Fig. 1* Scheme of the capsule drive system. $M$ - mass of the capsule, $m$ – mass of the pendulum, $l$ – length of the pendulum, $\theta$ – pendulum angle, $k$ – spring stiffness, $c$ – damping coefficient, $F_\theta$ – external torque acting on the pendulum, $F_x$ – friction force, $x(t), y(t)$ – coordinates of the capsule [15]

In the system under consideration, propulsion of the capsule is caused by interactions between friction $F_x$, horizontal inertia force produced by swinging motion of the pendulum $R_x$ and reaction (contact) force between the capsule and the underlying surface $R_y$. Dynamics of the presented system is thoroughly described in [10–12]. Moreover, detailed derivation of equations of motion of the pendulum capsule can be found in our previous work [16]. Therefore, in this paper, we only cite the final form of equations of motion. The following non-dimensional quantities are used:

$$\Omega = \sqrt{\frac{g}{l}}, \tau = \Omega t, \gamma = \frac{M}{m}, z = \frac{x}{l}, \rho = \frac{k}{m\Omega^2 l^2}, \nu = \frac{c}{m\Omega l^2},$$
$$f_z = \frac{F_x}{m\Omega^2 l}, u = \frac{F_\theta}{m\Omega^2 l^2}, r_z = \frac{R_x}{m\Omega^2 l}, r_y = \frac{R_y}{m\Omega^2 l}, \quad (1)$$

where $g$ is the gravitational acceleration, $t$ and $\tau$ correspond to dimensional and dimensionless time, respectively, and all the other symbols are described in the caption of Fig. 1. Relations between derivatives with respect to $t$ and $\tau$ are as follows.

$$\dot{x} = \frac{dx}{dt} = \frac{dx}{d\tau}\frac{d\tau}{dt} = \Omega\frac{dx}{d\tau} = \Omega x', \quad \ddot{x} = \frac{d^2x}{dt^2} = \frac{d}{dt}\left(\frac{dx}{dt}\right) = \frac{d}{d\tau}\left(\Omega\frac{dx}{d\tau}\right)\frac{d\tau}{dt} = \Omega^2\frac{d^2x}{d\tau^2} = \Omega^2 x''. \quad (2)$$

Using symbols and notation defined in the formulas (1), (2), equations of motion of the pendulum capsule drive can be presented in the following dimensionless, matrix form [16].

$$\begin{bmatrix} 1 & -\cos\theta(\tau) \\ -\cos\theta(\tau) & \gamma+1 \end{bmatrix}\begin{bmatrix} \theta''(\tau) \\ z''(\tau) \end{bmatrix} = \begin{bmatrix} \sin\theta(\tau) - \rho\theta(\tau) - \nu\theta'(\tau) + u(\tau) \\ -\theta'^2(\tau)\sin\theta(\tau) - f_z(\tau) \end{bmatrix}. \quad (3)$$

Further dimensionless quantities are contact force $r_y$, resultant horizontal load due to pendulum's motion $r_z$ and the dimensionless Coulomb friction $f_z$, that are described by the following equations:

$$r_y(\tau) = (\gamma + 1) - \theta''(\tau)\sin\theta(\tau) - \theta'^2(\tau)\cos\theta(\tau), \quad (4)$$

$$r_z(\tau) = \theta''(\tau)\cos\theta(\tau) - \theta'^2(\tau)\sin\theta(\tau), \quad (5)$$

$$f_z(\tau) = \begin{cases} \mu r_y(\tau)sgn[z'(\tau)] \leftrightarrow z'(\tau) \neq 0, \\ \mu r_y(\tau)sgn[r_z(\tau)] \leftrightarrow z(\tau) = 0 \wedge |r_z(\tau)| \geq \mu r_y(\tau), \\ r_z(\tau) \leftrightarrow z(\tau) = 0 \wedge |r_z(\tau)| < \mu r_y(\tau), \end{cases} \quad (6)$$

where $\mu$ is the friction coefficient. Equations (3)-(6) form the complete model of the capsule pendulum drive presented in Fig. 1.

Our previous paper [16] describes a numerical method which enables to approximate the optimal control of a system in the form of a finite number of Fourier series terms (7).

$$u(\tau) = \frac{a_0}{2} + \sum_{k=1}^{K} a_k \cos(k\omega\tau) + \sum_{k=1}^{K} b_k \sin(k\omega\tau). \quad (7)$$

Control of the system (3)-(6) has been optimized with respect to the distance covered by the capsule within the dimensionless time interval $\tau \in [0, 100]$. In optimization process, the following values of system parameters have been assumed: $\mu = 0.3, \rho = 2.5, \nu = 1.0, \gamma = 10$. Moreover, it has been asserted that the control $u(\tau)$ has to remain in the allowable range $[-4, 4]$. Under such assumptions, taking $K = 5$ harmonics in the formula (7), the following parameters of the approximate, open-loop optimal control have been obtained.

$$\begin{aligned} a_0 &\approx 1.62506, \quad \omega \approx 1.64722 \\ (a_1, a_2, a_3, a_4, a_5) &\approx (-3.43222, -1.95285, -0.68182, 0.38493, 0.17389) \\ (b_1, b_2, b_3, b_4, b_5) &\approx (-0.41690, 0.12411, -0.10468, 0.13722, 0.27902) \end{aligned} \quad (8)$$

The control function (7) with parameters (8) is plotted in Fig. 4 and the resulting trajectory of capsule's motion is presented in Fig. 5 in the following sections of the paper. These results are the starting point for the current research. In the remaining part of this work, we are going to show that a neural network can learn from the open-loop solution (7) in order to form a closed-loop controller. Moreover, it is going to be presented that, counterintuitively, such a neural network can outperform the original solution which served as the training set in the learning process (see Fig. 5).

## 3. Research methodology

In this chapter, the use of a neural network as a closed-loop controller based on the approximate solution of the open-loop optimal control is described. The optimized open-loop control function is obtained by means of the Fourier series based method [16].

The research is divided into three stages. The first one, named preliminary research results, focuses on the evaluation of the learning process in which the neural network approximates dependencies between the optimized open-loop control and corresponding states of the system. To achieve this, a neural network model consisting of the three layers was created in Python language and tested for various hidden and output layer activation functions, along with the changing number of neurons. One chosen optimizer was applied to all tests. For each model of the neural network, learning process has been repeated 3 times in order to reduce the risk of randomness of results. Neural network model performance was assessed with regard to the data fit and the loss function value with the use of the two chosen parameters described in the following chapters.

Results with the top score of model performance were used in the second stage. In this part, the obtained neural network parameters have been implemented in the simulated controller of the capsule drive. This allowed to calculate the distance covered by the capsule pendulum system driven by the neural controller.

In the last part, robustness of both controllers (the optimized open-loop and the neural closed-loop ones) with respect to variation of the friction coefficient between the capsule and the underlying surface is examined. Performance of the controllers is tested in different conditions, from constant friction coefficient to large variations of this parameter.

### 3.1. Neural Network structure

Consider a discontinuous dynamical system, for which the optimal trajectory under a Fourier series-based numerical method of open-loop control optimization is calculated [16]. The objective of the neural network is to return the value of optimized control for an arbitrary state of the controlled object (i.e. the pendulum capsule drive). Therefore, input values of the network are state variables of the system (3)-(6) and the output value is expected to be the corresponding control.

The neural network structure design demands few key decisions to be made, such as number of layers, number of neurons and activation functions. The number of hidden layers depends on the problem with which it deals. In this case, one hidden layer has been fixed, due to the fact that it is enough to approximate any arbitrary continuous function [26–27].

Another problem to be solved is the number of neurons. It should be chosen in a way avoiding the under or overfitting. In the first case, where there are not enough neurons to train, results could not be satisfying and adequate to what is expected. On the other hand, with too many neurons, the artificial neural network possesses an exceeding number of parameters to be determined, which makes it "remember" each data point separately and lose the data generalization property. There exist many rules of thumb for establishing the correct number of neurons in a hidden layer [26]. The authors decided to test one of them, where the number of hidden neurons is calculated as the sum of the 2/3 size of the input layer and the whole size of the output layer. In total, five different numbers of neurons were tested, with four values chosen randomly.

One of the exemplary activation functions used for the hidden layer is Rectified Linear Unit (ReLU) which is considered one of the most efficient due to its good resistance to the vanishing gradients [28]. In this research, the authors used also the sigmoid (logistic) and hyperbolic tangent functions. The last one is preferred more since its gradients are not restricted to vary in the specified direction [29]. Moreover, when the output returned by the sigmoid function is close to zero caused by the highly negative inputs, the process of the neural network elongates and the probability of getting stacked in some local minima is higher [30].

For the output layer, a default option is a linear function, commonly used for regression problems. However, authors tested sigmoid function in the output layer too.

Since the structure of the neural network model consisting of the specified number of neurons and hidden layers, as well as the activation functions for each layer, is established, the default neural network model for *n* neurons proposed by the authors is presented in Fig. 2.

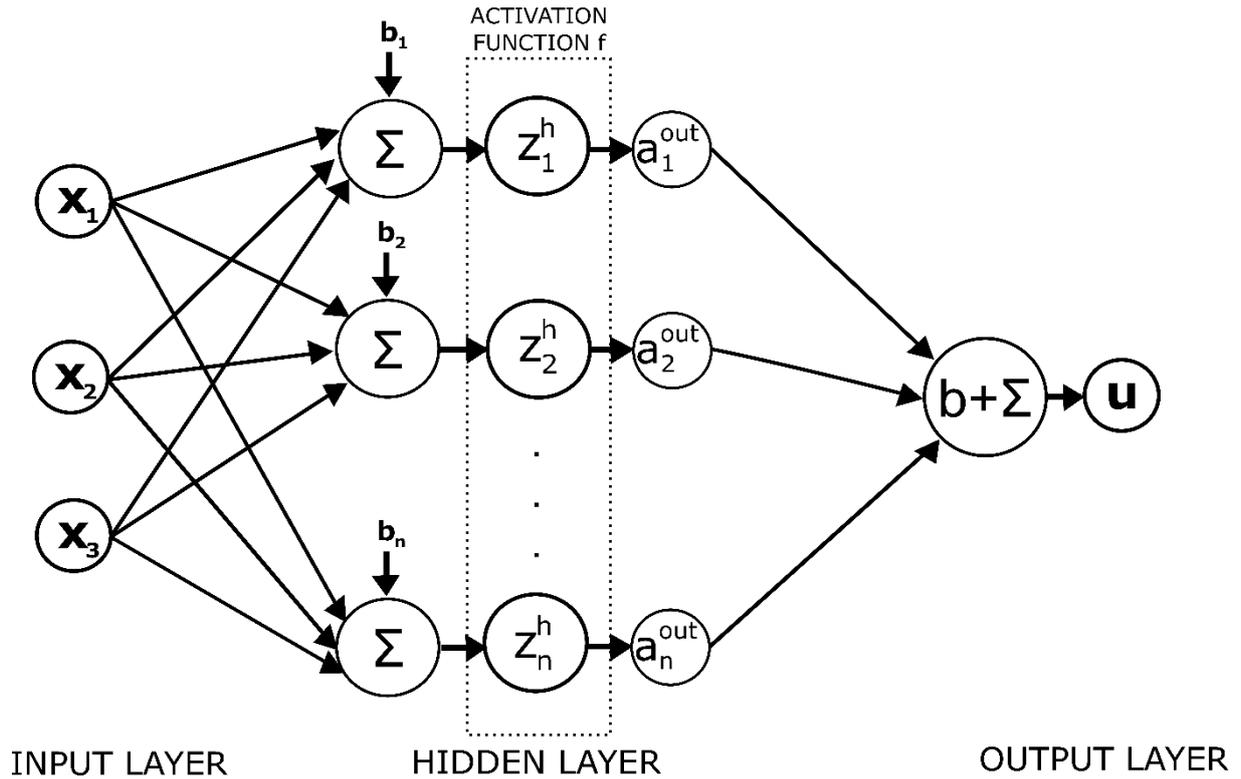

**Fig. 2** Neural network structure, where: $n$ – number of current neuron, $x_1$, $x_2$, $x_3$ – system state variables and input activation for the hidden layer, $b, b_1, b_2, \ldots b_n$ – biases for the current neuron in the hidden layer, $z_n^h$ – value for the current neuron in the hidden layer, $a_1^{out}, a_2^{out}, \ldots a_2^{out}$ – input activation for the output layer, $u$ – output (control)

The proposed neural network consists of 3 layers. An input layer is represented by a vector $x$ with shape (3,1). In this case, $x_1$ refers to the position of the pendulum $\theta$, whereas $x_2$ and $x_3$ represents the velocity of the pendulum $\dot{\theta}$ and capsule system $z'$, respectively. Please note that position of the capsule $z$ does not appear here, as the control is expected to be independent from this variable. In the (unique) hidden layer, the activation function $f$ is applied. It could be imagined as a filter that the values going through scaling the output into the proper ranges. Without any applied activation function, the neural network could learn only linear transformation. The output layer consists of one neuron, which represents the capsule pendulum drive system trajectory $u$ that will be predicted by the neural network. For each neuron in the hidden layer, the input value $a_n^{out}$ for the last layer is calculated in the following manner:

$$a_n^{out} = f(z_n^h) = f\left(b_n^h + \sum_{i=1}^{P} x_i w_{ni}\right) \tag{9}$$

where:

$n$ – index of the current neuron,

$h$ - hidden layer,

$f$ – the activation function applied in the hidden layer (e.g. ReLU, sigmoid, tanh),

$z_n^h$ - value for the current neuron in the hidden layer,

$b_n^h$ – bias for the current neuron in the hidden layer,

$x_i, a_i$ – an input activation for the following layer, i.e. an input activation for the hidden layer is represented by three system state variables $x_1, x_2$ and $x_3$,

$w_{ni}$ – weight representing the connection between the current neuron and the neurons from the previous layer,

$i$ – the number of the current neuron from the previous layer,

$P$ - number of inputs activation/neurons from the previous layer.

## 3.2. Learning method and tools

In the preliminary research, the script for the neural network learning process in Python language was created. Firstly, the optimal trajectory of the capsule pendulum drive system obtained with the use of the Fourier series based optimization algorithm that is the aim of the neural network prediction was loaded and checked concerning missing values. Uploaded input (system state variables) and output (trajectory) data were divided into the training and test sets in the 80:20% proportion providing necessary information for the learning process. In this step, the authors implemented two additional parameters. The first one was responsible for dataset shuffling, which provides to obtain representative training and test sets. It means also that the created model is not just as it is because of the applied data order. Moreover, it allows preventing getting stuck in cycles during the cost function optimization [30]. The another parameter was used to determine the seed for dataset shuffling, guaranteeing consistency and reproducible results [30]. The default neural network model proposed by the authors, visible in Fig. 2 was implemented in the script.

The way of the how the loss gradients is used to update the parameters of the neural network is specified by the optimizers [23]. In this research, the authors consider the adaptive moment estimation (Adam) optimizer combining the advantages of AdaGrad [31] and RMSProp [32] methods [33]. The first one deals efficiently with sparse gradients, whereas the second works well in online, as well as non-stationary settings and resolves some problems of the first one. More precise connections between these methods and Adam are described in [33]. In the chosen optimizer, the hyper-parameters are equipped with intuitive interpretation and typically do not require any tuning. The individual adaptive learning rates are calculated for different parameters based on the estimation of the first and second moments of the gradient [33].

When the optimizer is implemented, the process of data prediction can be performed. To avoid overtraining, the loss values reaching the same level were monitored for the specified value and stopped consequently. The process of neural network learning was fixed for 1000 epochs. The performance of the model was assessed with regard to the loss function value and data fit on the test dataset. The value of the loss function was measured with the use of the mean squared error (MSE), which is the average value of the cost function squared errors that is minimized to fit the model [30]. Data fit evaluation how far the predicted values from the original one are, was described using the coefficient of determination $R^2$ score meaning the fraction of response variance captured by the model [30].

## 4. Research results

The evaluation of the neural network learning the optimal Fourier series based trajectory process has been analyzed for 30 artificial models being the modification of the default one proposed by the authors, with various activation functions for the hidden and output layers, as well as different numbers of neurons, with the use of the Adam optimization algorithm. Thirty neural network models were compared regards to the data fit and loss function score represented by $R^2$ and MSE, respectively. The top values describing the neural network performance, obtained for various sets of artificial models are presented in Table 1.

**Table 1** *The top score of a data fit ($R^2$) and loss function (MSE) refer to the predicted data for the thirty various sets of neural network model optimized by the Adam algorithm, consisting of the changing number of neurons and three layers, where different activation functions for the hidden and output layer have been applied*

| OPTIMIZER ADAM | Output layer activation function | None (Linear) | | Sigmoid | |
|---|---|---|---|---|---|
| Hidden layer activation function | | Data fit/Loss function score | | | |
| | No of neurons | $R^2$ | MSE | $R^2$ | MSE |
| ReLU | 3 | 0.985 | 0.002 | 0.985 | 0.0019 |
| | 10 | 0.996 | 0.0005 | 0.995 | 0.0006 |
| | 17 | 0.997 | 0.0003 | 0.998 | 0.0003 |
| | 30 | 0.998 | 0.0002 | 0.998 | 0.0003 |
| | **50** | **0.999** | **0.0001** | **0.999** | **0.0002** |
| Sigmoid | 3 | 0.991 | 0.0012 | 0.989 | 0.0014 |
| | 10 | 0.996 | 0.0005 | 0.998 | 0.0003 |
| | **17** | 0.997 | 0.0004 | **0.999** | **0.0001** |
| | **30** | 0.997 | 0.0005 | **0.999** | **0.0001** |
| | **50** | 0.998 | 0.0003 | **0.999** | **0.0001** |
| Tanh (hyperbolic tangent) | 3 | 0.965 | 0.0045 | 0.989 | 0.0015 |
| | **10** | 0.996 | 0.0005 | **0.999** | **0.0001** |
| | 17 | 0.998 | 0.0003 | **0.999** | **0.0002** |
| | **30** | 0.998 | 0.0003 | **0.999** | **0.0001** |
| | **50** | 0.996 | 0.0006 | **0.999** | **0.0001** |

From the preliminary research stage, the authors decided to choose the neural network models with the performance score at the level of 0.999 and 0.0001/0.0002, for the $R^2$ and MSE, respectively. Consequently, nine different sets (bolded in Table 1) with different numbers of neurons and activation functions applied for hidden and output layers reached this result. It can be seen that results with the highest scores are mainly obtained for the neural network models with sigmoid and hyperbolic tangent activation functions for the hidden layer, whereas for the output one the sigmoid is applied. Moreover, high-scored performance results with the mentioned combination of activation functions for hidden and output layers start at 10 neurons without significant changes during the further increase in value. On the other hand, the thumb rule tested by authors in this study does not give satisfying results, which could be related to the specific character of the considered system with the small number of input and output data.

The linear relationship between the learning and data predicted from the neural network learning process, along with the $R^2$ coefficient value equal to 0.999, showing the high data fit quality, is presented in Fig. 3. Results were obtained for the neural network model consisting of 50 neurons in the hidden layer along with the ReLU applied as an activation function and a linear one for the output.

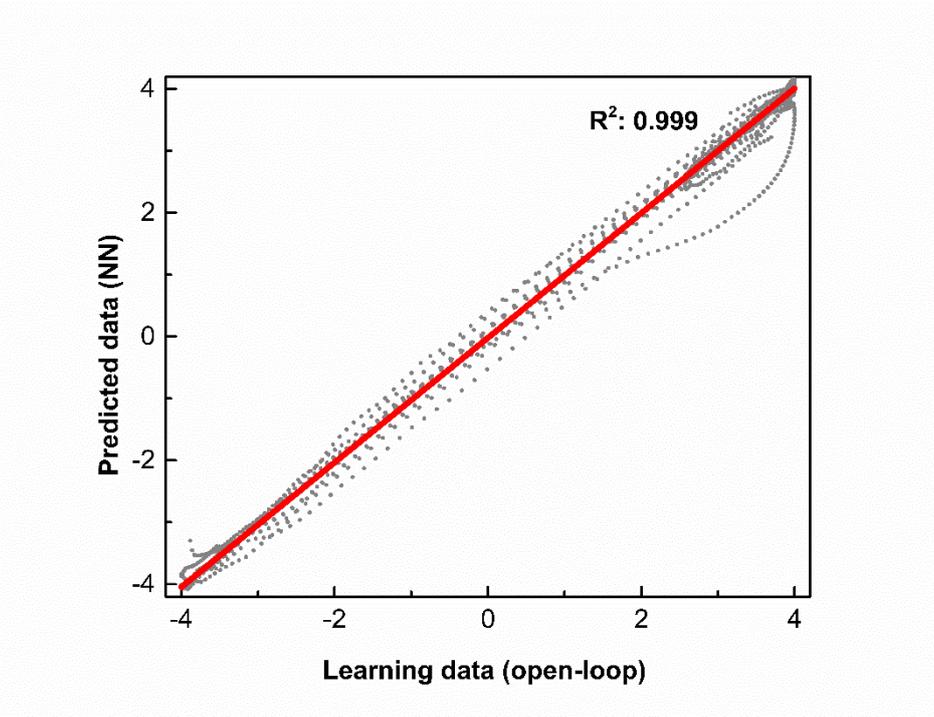

*Fig. 3* *Linear relationship between the learning and predicted data along with the $R^2$ coefficient value equal to 0.999 for the neural network model consisting of 50 neurons, ReLU applied as an activation function for the hidden layer and linear one for the output, respectively*

Parameters of the neural network models with the highest performance obtained in the preliminary research results have been implemented into the capsule pendulum drive controller to calculate the distance covered by the system in the dimensionless time interval. Results have been analyzed and presented in Table 2.

*Table 2* Distance covered by the capsule pendulum drive system for the neural network models with the highest performance score obtained in the preliminary research along with their structure including the hidden and output layer activation functions, number of neurons, the data fit and loss function scores

| Hidden layer's activation function | Output layer's activation function | No. of neurons | Data fit score | Loss function score | Distance |
|---|---|---|---|---|---|
| ReLU | None (linear) | 50 | 0.999 | 0.0001 | 6.135 |
| ReLU | Sigmoid | 50 | 0.999 | 0.0002 | 5.998 |
| Sigmoid | Sigmoid | 17 | 0.999 | 0.0001 | 6.081 |
| Sigmoid | Sigmoid | 30 | 0.999 | 0.0001 | 6.082 |
| Sigmoid | Sigmoid | 50 | 0.999 | 0.0001 | 6.073 |
| Tanh | Sigmoid | 10 | 0.999 | 0.0001 | 6.115 |
| Tanh | Sigmoid | 17 | 0.999 | 0.0002 | 6.079 |
| Tanh | Sigmoid | 30 | 0.999 | 0.0001 | 6.076 |
| Tanh | Sigmoid | 50 | 0.999 | 0.0001 | 6.089 |

The analysis of the distance covered by the capsule pendulum drive system in the dimensionless time interval, tested for the various sets of neural network models, revealed that the highest score is obtained for the neural network consisting of 50 neurons in the hidden layer along with the ReLU applied as an activation function and linear one used for the output. The achieved result is equal to 6.135 with 1.16% higher performance than the one from the open-loop control (6.065). The other neural network models scores in this stage differ slightly from the reference open-loop trajectory control and have not been considered in this study. It is in any case worthy to observe that all (but one) provide a final displacement that is systematically larger of the reference one 6.065, although the increment is minor.

The comparison of the control trajectory, as well as the distance covered by the capsule pendulum drive system in the dimensionless time interval obtained for both controllers (open-loop and neural closed-loop) is presented in Fig. 4 and Fig. 5, respectively.

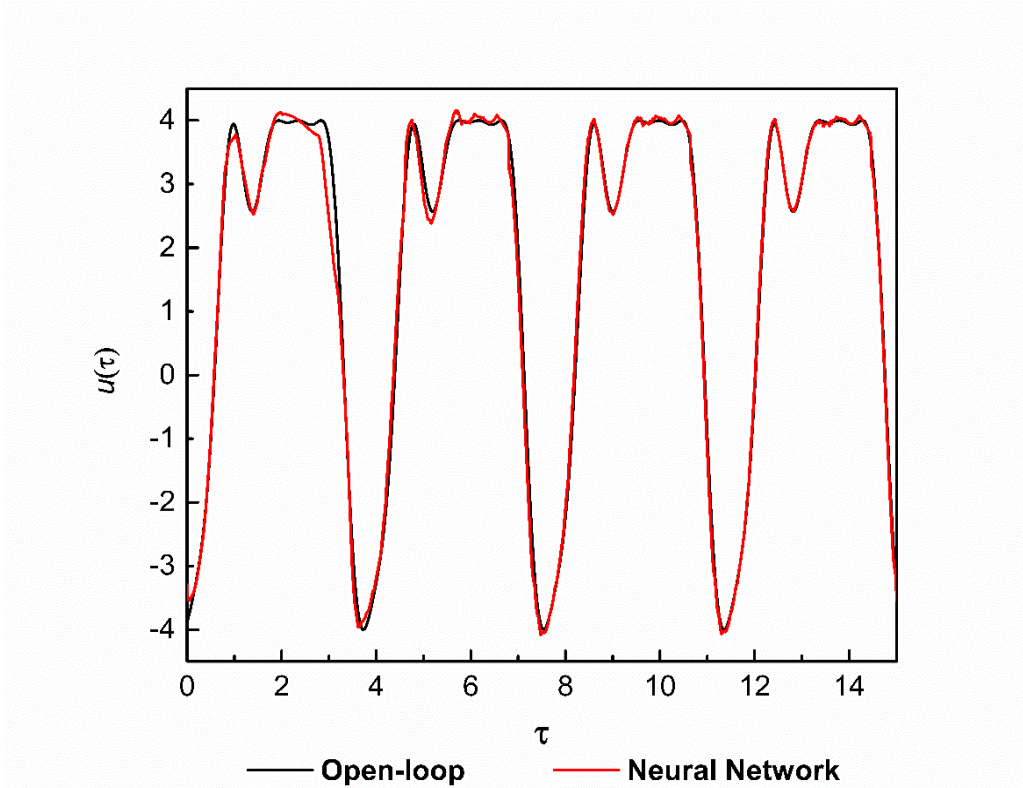

*Fig. 4* Control trajectory for the open-loop and neural network controller vs. dimensionless time ($\tau$)

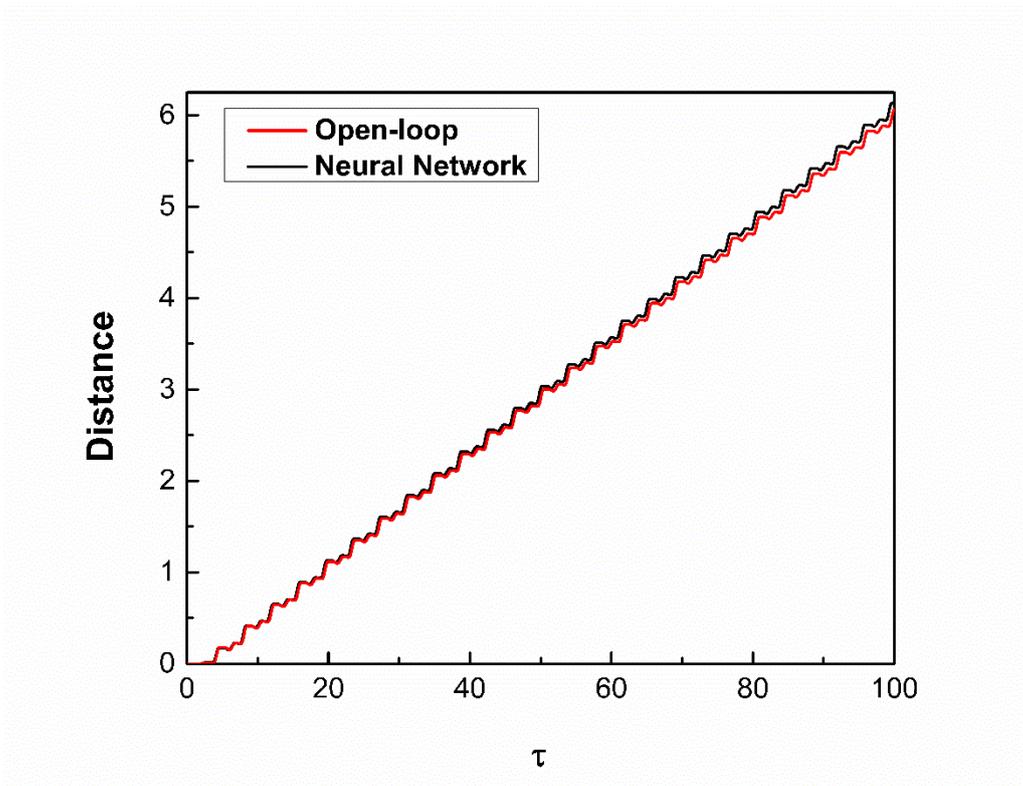

*Fig. 5* Distance covered by the capsule pendulum drive system for the open-loop and neural network controller vs. dimensionless time ($\tau$)

## 5. Perturbations in the system

The robustness of controllers against variation of the system's parameters has been tested by introducing a variable friction coefficient between the capsule and the underlying surface. It has been assumed that the actual friction coefficient at each point is a uniformly distributed random value:

$$\mu_r \in [\mu - \Delta, \ \mu + \Delta] \tag{10}$$

where $\mu$ is the nominal, assumed value of the friction coefficient and $\Delta$ is its maximal, absolute deviation. A different value of $\mu_r$ has been drawn for each interval of the path of motion of the capsule of the dimensionless length equal to 0.1.

The subject of the evaluation was the neural-network model with the highest obtained result of the distance covered by the capsule pendulum drive system in the interval of dimensionless time. The $\Delta$ parameter introducing perturbations was varying from 0.00 to 0.20, changing the effective friction coefficient range simultaneously. The mean score for each value of the $\Delta$ parameter with the corresponding standard deviation is presented in Table 3. Moreover, the relative changes between non- and perturbed distance for the open-loop and neural closed-loop controllers, as well as their comparison in the uncertain frictional environment are presented too.

*Table 3* Open-loop and neural network controllers comparison due to the mean distance covered by the capsule pendulum drive system in the uncertain frictional environment represented by the varying Δ parameter in the 0.00-0.20 range with the corresponding standard deviation. Relative changes between the non- and perturbed distance for both controllers along with the relative changes across both controllers working in the environment with the varying friction coefficient range. Parameters of the open-loop controller were calculated regards to the formula (7) and presented in (8). The neural network controller parameters obtained for the structure consisting of 50 neurons in the hidden layer along with the ReLU applied as an activation function and linear one for the output were calculated according to the formula (9)

| | Open-loop controller | | Closed-loop Neural Network controller | | Neural Network and open-loop comparison |
|---|---|---|---|---|---|
| | Distance without perturbation | | | | |
| | 6.065 | | 6.135 | | |
| $\Delta$ | Mean (±SD) distance affected by perturbations | Relative changes between non- and perturbed distance [%] | Mean (±SD) distance affected by perturbations | Relative changes between non- and perturbed distance [%] | Relative changes across both controllers' results [%] |
| 0.00 | 6.065 (± 0) | 0 | 6.135 (± 0) | 0 | 1.16 |
| 0.01 | 6.069 (± 0.0094) | 0.07 | 6.139 (± 0.0026) | 0.06 | 1.16 |
| 0.02 | 6.074 (± 0.012) | 0.15 | 6.138 (± 0.0044) | 0.03 | 1.04 |
| 0.03 | 6.040 (± 0.024) | -0.41 | 6.133 (± 0.040) | -0.04 | 1.54 |
| 0.04 | 6.068 (± 0.018) | 0.06 | 6.159 (± 0.038) | 0.39 | 1.49 |
| 0.05 | 6.087 (± 0.056) | 0.37 | 6.130 (± 0.015) | -0.09 | 0.70 |
| 0.06 | 6.042 (± 0.034) | -0.38 | 6.108 (± 0.040) | -0.45 | 1.10 |
| 0.07 | 5.999 (± 0.031) | -1.08 | 6.139 (± 0.033) | 0.06 | 2.33 |
| 0.08 | 6.048 (± 0.019) | -0.27 | 6.087 (± 0.057) | -0.79 | 0.64 |
| 0.09 | 5.978 (± 0.037) | -1.44 | 6.079 (± 0.029) | -0.93 | 1.69 |
| 0.10 | 5.976 (± 0.076) | -1.47 | 6.076 (± 0.039) | -0.98 | 1.67 |
| 0.11 | 5.951 (± 0.043) | -1.88 | 6.162(± 0.086) | 0.44 | 3.55 |
| 0.12 | 5.904(± 0.054) | -2.66 | 6.151 (± 0.32) | 0.25 | 4.19 |
| 0.13 | 5.846 (± 0.050) | -3.61 | 6.214 (± 0.18) | 1.28 | 6.29 |
| 0.14 | 5.873 (± 0.11) | -3.16 | 6.114 (± 0.085) | -0.35 | 4.10 |
| 0.15 | 5.864 (± 0.062) | -3.32 | 6.058 (± 0.11) | -1.27 | 3.31 |
| 0.16 | 5.775 (± 0.098) | -4.78 | 6.179 (± 0.068) | 0.71 | 6.99 |
| 0.17 | 5.803 (± 0.11) | -4.32 | 6.020 (± 0.38) | -1.88 | 3.74 |
| 0.18 | 5.638 (± 0.044) | -7.05 | 5.929 (± 0.078) | -3.36 | 5.17 |
| 0.19 | 5.614 (± 0.071) | -7.44 | 5.736 (± 0.14) | -6.50 | 2.19 |
| 0.20 | 5.443 (± 0.251) | -10.25 | 5.584 (± 0.12) | -9.00 | 2.57 |

The capsule pendulum drive system distance covered in the uncertain frictional environment where the $\Delta$ parameter was introduced to change the friction coefficient range between the capsule and underlying surface, presented in Table 3, shows that the neural network controller is more resistant to the occurred friction changes than the open-loop controller. The first significant decrease in the distance appears much earlier in the open-loop control, starting from the value equal to 0.07, corresponding to the 0.10 - 0.13 in the neural network controller. As the $\Delta$ parameter increase, the difference between the distances obtained for both controllers is more noticeable. The close-up look of the distance consideration for open-loop and neural network controllers is presented in Fig. 6.

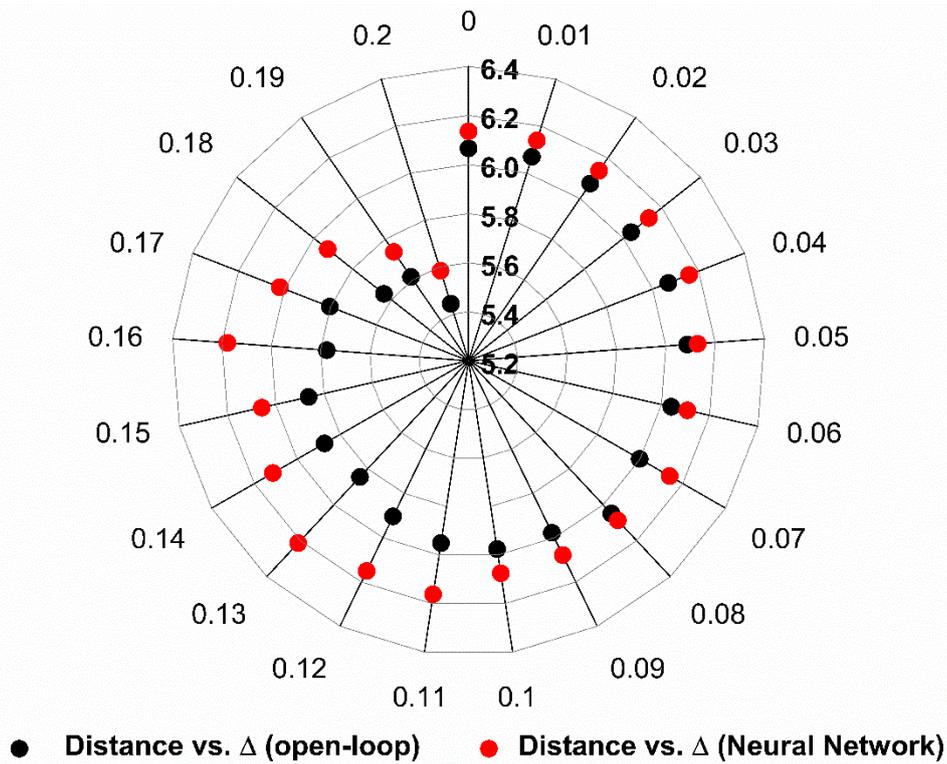

*Fig. 6 Distance covered by the capsule pendulum drive system in the uncertain frictional environment for the open-loop and neural network controller*

To better understand the changes taking place in the system working in the environment with varying friction coefficient range, the close-up look at the relative changes between the non- and perturbed distances, as well as the comparison of both controllers regards to the perturbed distances is shown too (see Fig. 7). It can be noticed that the higher the $\Delta$ value, the more decrease in the distance is observed, especially for the open-loop control. Moreover, the level of the 1% percent changes between the non- and perturbed distance for the neural network controller remains till the $\Delta$ parameter value is equal to 0.16 with two jumps before, whereas the changes for the open-loop controller pass this level at least 4 times simultaneously. The greatest change observed for both controllers appears when the $\Delta$ parameter is equal to 0.20, meaning the 10.25% and 9% decrease of the distance for the open-loop controller and the neural network one, respectively. While for large values of $\Delta$ the performances of both controllers decrease (although to a different extent, as said), it is interesting to note that, quite surprisingly, for very small values of $\Delta$ the maximum distance increases with respect to the unperturbed case $\Delta = 0$, suggesting a kind of beneficial effect of uncertainties on $\mu$.

The perturbed distances comparison analysis performed for both controllers proves that the neural network controller is more resistant to the changes occurring in the environment with varying friction coefficient range. In this case, the most significant observed increase of performance is equal to +7%.

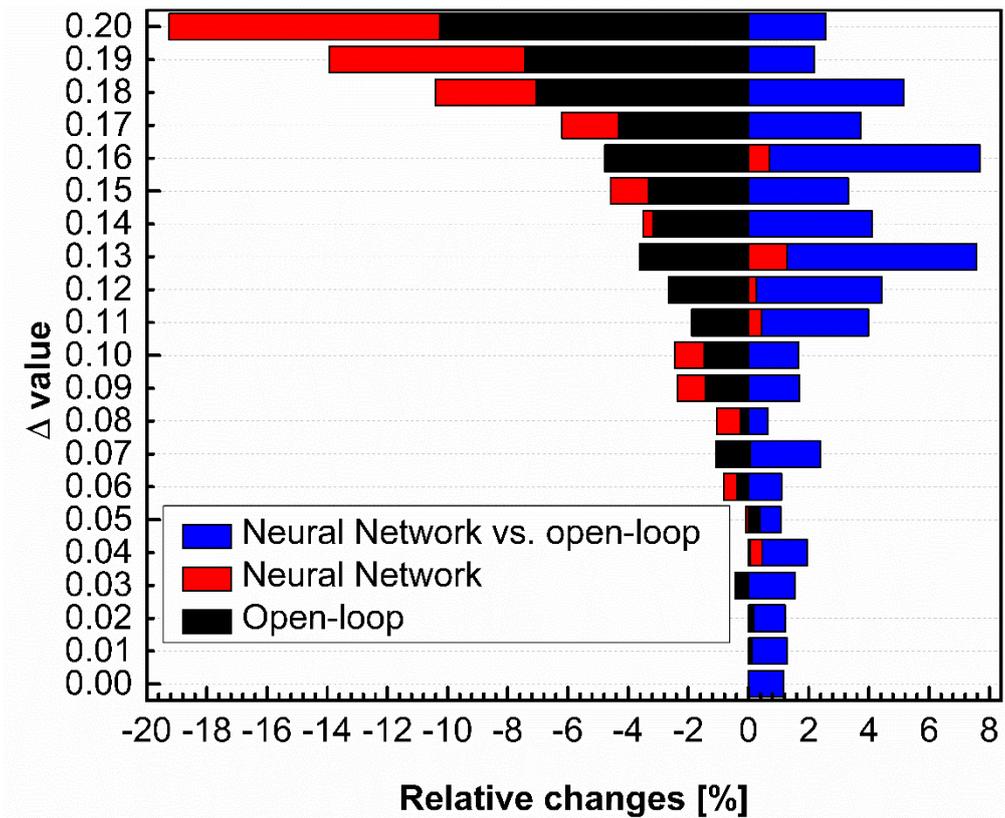

*Fig. 7* Relative changes between non- and perturbed distances for the open-loop and neural network controllers with their comparison regards to the work in the environment with varying friction coefficient range introduced by Δ parameter value

## 6. Summary and conclusions

In this study, a novel approach of the pendulum capsule drive system control with the use of the neural network as a closed-loop controller is presented. The novelty in this research is the optimized open-loop control function being the base of the designed controller. The dependence between the output of the open-loop controller and the system state variables is determined by the neural network.

The main aim of this research was to test and evaluate the robustness of the novel controller compared with the original open-loop control one. Thus, the study was divided into three parts. In the preliminary research, the authors created a model of a neural network that was tested for various sets with different numbers of neurons, as well as the activation functions for hidden and output layers. In this stage, the performance of each model measured with the $R^2$ and MSE tools allowed to select nine neural networks models that have the highest correlation between the learning and predicted data. Gathered parameters of the neural network learning process have been implemented in the capsule pendulum drive simulated controller to calculate the distance covered by the system in the dimensionless time interval. The highest achieved distance was obtained for the neural network model consisted of 50 neurons in the hidden layer along with the ReLU used as an activation function and a linear one applied for the output.

This artificial model was tested in the last part of the study considering the robustness of the controller and its reliability in the uncertain frictional environment, introduced by the varying friction coefficient range between the capsule shell and the underlying surface. With the constant friction coefficient, the neural network controller revealed a 1.16% higher performance than the open-loop. Further tests performed in the expanded friction range coefficient variation proved that the neural network is more resistant to the perturbations occurring in the system, with a maximum of +7% advantage over the open-loop controller. Additionally, the changes between the non- and perturbed distance occurred there much slower and remained with the 1% change much longer. Meanwhile, in the open-loop controller, the increasing value was constantly observed exceeding some trial scores achieved in the neural network controller at least 4 times.

Results presented in this study confirm that the neural controller works more efficiently compared to the original open-loop controller and proves the higher level of robustness in the environment where perturbations occur. It seems that the neural network controller could be an alternative option to the classic ones in many applications of the mechanical field, especially for non-smooth and discontinuous systems (as the one considered in this work). Moreover, it could simplify significantly the closed-loop controllers' systems design where the open-loop control is only available.

## Acknowledgements

This study has been supported by the National Science Centre, Poland, PRELUDIUM Programme (Project No. 2020/37/N/ST8/03448).

This study has been supported by the National Science Centre, Poland under project No. 2017/27/B/ST8/01619.

This paper has been completed while the first author was the Doctoral Candidate in the Interdisciplinary Doctoral School at the Lodz University of Technology, Poland.

Sandra Zarychta, Marek Balcerzak and prof. Andrzej Stefański declare that their research described in this work has been carried out during their scientific visit at Marche Polytechnic University, Ancona, in collaboration with prof. Stefano Lenci.

The majority of the work has been done during the first author's three months stay at Marche Polytechnic University, Ancona.

## Conflict of interest

The authors declare that they have no conflict of interest.

## Data availability

Data obtained in this study are available from the corresponding author on reasonable request.